\documentclass[11pt]{amsart}

\usepackage{amsmath,amsthm,amssymb,url}

\usepackage{fancyhdr}

\usepackage[all]{xy}

\DeclareMathOperator{\Hom}{Hom}

\DeclareMathOperator{\Def}{Def}
\DeclareMathOperator{\End}{End}

\DeclareMathOperator{\MC}{MC}

\DeclareMathOperator{\holim}{holim}

\newcommand{\T}{\mathcal{T}}

\newcommand{\Eps}{\mathcal{E}}
\newcommand{\lie}{\boldsymbol{l}}

\newtheorem*{theorem*}{Theorem}
\newtheorem*{proposition*}{Proposition}
\newtheorem*{corollary*}{Corollary}

\theoremstyle{definition}

\begin{document}

\title[$\infty$-groupoids and period maps]{A short note on $\infty$-groupoids and the period map for projective manifolds}


\author{Domenico Fiorenza}

\address{Dipartimento di Matematica ``Guido
Castelnuovo''\\
Universit\`a di Roma ``La Sapienza''\\
P.le Aldo Moro 5,
I-00185 Roma Italy.}
\email{fiorenza@mat.uniroma1.it}

\author{Elena Martinengo}
\address{Institut f\"ur Mathematik und Informatik\\
Freie Universit\"at Berlin\\
         Arnimallee 3\\
         14195 Berlin, Germany}
\email{elenamartinengo@gmail.com}

{\let\thefootnote\relax\footnotetext{Received July 30th, 2011; revised August 18th, 2012; published September 13th, 2012.}}

\maketitle

\thispagestyle{fancy}

\begin{abstract}
We show how several classical results on the infinitesimal behaviour of the period map for smooth projective manifolds can be read in a natural and unified way within the framework of $\infty$-categories.
\end{abstract}

\section*{}
A common criticism of $\infty$-categories in algebraic geometry is that
they are an extremely technical subject, so abstract to be useless in everyday
mathematics. The aim of this note
is to show in a classical example that quite the converse is true: even a na\"{\i}ve
intuition of what an $\infty$-groupoid should be clarifies several aspects of the
infinitesimal behaviour of the periods map of a projective manifold. In
particular, the notion of Cartan homotopy turns out to be completely natural from this perspective, and so classical results such as Griffiths' expression for the differential
of the periods map, the Kodaira principle on obstructions to
deformations of projective manifolds, the Bogomolov-Tian-Todorov theorem, and Goldman-Millson quasi-abelianity theorem are easily recovered. 
\par
The use of the language of $\infty$-categories should not be looked at as providing new proofs for these results; namely, up to a change in language, our proofs verbatim reproduce arguments from the recent literature on the subject, particularly from the work of Marco Manetti and collaborators on dglas in deformation theory. Rather, by this change of language we change our point of view on the classical theorems above: in the perspective of $\infty$-sheaves from \cite{lurie}, all these theorems have a very simple \emph{local} nature which can be naturally expressed in terms of $\infty$-groupoids (or, equivalently, of dglas); their classical \emph{global} counterparts are then obtained by taking derived global sections. It is worth remarking that, if one prefers proofs which do not rely on the abstract machinery of $\infty$-categories, one can rework the arguments of this note in purely classical terms. Namely, once the abstract $\infty$-nonsense has suggested the ``correct'' local dglas, one can globalize them by means of an explicit model for the derived global sections, e.g., via resolutions by fine sheaves as in \cite{fiorenza-manetti2}, or by the Thom-Sullivan-Whitney model as in \cite{iacono-manetti}. 

 \par
Since most of the statements and constructions we recall in the paper are well known in the $(\infty,1)$-categorical folklore, despite our efforts in giving credit, it is not unlikely we may have misattributed a few of the results; we sincerely apologize for this. We thank the referee for accurate  remarks which helped us a lot in improving the present paper, and Ezra Getzler, Donatella Iacono, Marco Manetti, Jonathan Pridham, Carlos Simpson, Jim Stasheff, Bruno Vallette, Gabriele Vezzosi, and the $n$Lab for suggestions and several inspiring conversations on the subject of this paper. 
\par
Through the whole paper, ${\mathbb K}$ is a fixed characteristic zero field, all algebras are defined over ${\mathbb K}$ and local algebras have ${\mathbb K}$ as residue field. In order to keep our account readable, we will gloss over many details, particularly where the use of higher category theory is required.

\section{From dglas to $\infty$-groupoids and back again}

With any nilpotent dgla ${\mathfrak g}$ is naturally associated the simplicial set
\[
\MC({\mathfrak g}\otimes \Omega_\bullet),
\]
where $\MC$ stands for the Maurer-Cartan functor mapping a dgla to the set of its Maurer-Cartan elements, and $\Omega_\bullet$ is the simplicial differential graded  commutative associative algebra of polynomial differential forms on algebraic $n$-simplexes, for $n\geq0$. The importance of this construction, which can be dated back to Sullivan's \cite{sullivan}, relies on the fact that, as shown by Hinich and Getzler \cite{getzler, hinich}, the simplicial set $\MC({\mathfrak g}\otimes \Omega_\bullet)$ is a Kan complex, or -to use a more evocative name- an $\infty$-groupoid. A convenient way to think of $\infty$-groupoids is as homotopy types of topological spaces; namely, it is well known\footnote{At least in higher categories folklore} that any $\infty$-groupoid can be realized as the $\infty$-Poincar\'e groupoid, i.e., as the simplicial set of singular simplices, of a topological space, unique up to weak equivalence. Therefore, the reader who prefers to can substitute homotopy types of topological spaces for equivalence classes of $\infty$-groupoids. To stress this point of
view, we'll denote the $k$-truncation of an $\infty$-groupoid $\mathbf{X}$
by the symbol $\pi_{\leq k}\mathbf{X}$. More explicitely, $\pi_{\leq
k}\mathbf{X}$ is the $k$-groupoid whose $j$-morphisms are the $j$-morphisms of $\mathbf{X}$ for $j<k$, and are homotopy classes of $j$-morphisms of $\mathbf{X}$ for $j=k$. In
particular, if $\mathbf{X}$ is the $\infty$-Poincar\'e groupoid of a
topological space $X$, then $\pi_{\leq 0}\mathbf{X}$ is the set $\pi_0(X)$
of path-connected components of $X$, and $\pi_{\leq 1}\mathbf{X}$ is the
usual Poincar\'e groupoid of $X$. 
\par The next step is to
consider an $(\infty,1)$-category, i.e., an $\infty$-category whose
hom-spaces are $\infty$-groupoids. This can be thought as a formalization of
the na\"ive idea of having objects, morphisms, homotopies between morphisms,
homotopies between homotopies, et cetera. In this sense, endowing a category
with a model structure should be thought as a first step towards defining an
$(\infty,1)$-category structure on it.
\par
Turning back to dglas, an easy way to produce nilpotent dglas is the following: pick an arbitrary dgla ${\mathfrak g}$; then, for any differential graded local Artin algebra $A$, take the tensor product ${\mathfrak g}\otimes {\mathfrak m}_A$, where ${\mathfrak m}_A$ is the maximal ideal of $A$. Since both constructions 
\begin{align*}
{\bf DGLA}\times {\bf dgArt}&\to {\bf nilpotent\,\, DGLA}\\
(\mathfrak{g},A)&\mapsto \mathfrak{g}\otimes \mathfrak{m}_A
\end{align*}
and
\begin{align*}
{\bf nilpotent\,\, DGLA}&\to {\bf \infty\text{\bf -Grpd}}\\
\mathfrak{g}&\mapsto \MC({\mathfrak g}\otimes \Omega_\bullet)
\end{align*}
are functorial, their composition defines a functor
\[ \Def: {\bf DGLA} \to {\text{\bf Formal }\bf \infty\text{\bf -Grpd}},\]
where, by definition, a formal $\infty$-groupoid is a functor ${\bf dgArt}\to {\bf \infty\text{\bf -Grpd}}$. Note that $\pi_{\leq 0}(\Def({\mathfrak g}))$ is the usual set valued deformation functor associated with ${\mathfrak g}$, i.e., the functor
\[
A\mapsto \MC({\mathfrak g}\otimes {\mathfrak m}_A)\bigl/{\rm gauge},
\] 
where the gauge equivalence of Maurer-Cartan elements is induced by the gauge action
\[
e^\alpha*x=x+\sum_{n=0}^\infty \frac{({\rm ad}_{\alpha})^n}{(n+1)!}\ ([\alpha,x]-d\alpha)
\]
 of $\exp(\mathfrak{g}^0\otimes\mathfrak{m}_A)$ on the subset $\MC({\mathfrak g}\otimes {\mathfrak m}_A)$ of $\mathfrak{g}^1\otimes\mathfrak{m}_A$.
However, due to the presence of nontrivial irrelevant stabilizers, the groupoid $\pi_{\leq 1}(\Def({\mathfrak g}))$ is not equivalent to the action groupoid $\MC({\mathfrak g}\otimes {\mathfrak m}_A)\bigl/\bigl/\exp(\mathfrak{g}^0\otimes\mathfrak{m}_A)$, unless ${\mathfrak g}$ is concentrated in nonnegative degrees. We will come back to this later. Also note that the zero in $\mathfrak{g}^1\otimes\mathfrak{m}_A$ gives a natural distinguished element in $\pi_{\leq 0}(\Def({\mathfrak g}))$: the isomorphism class of the trivial deformation. Since this marking is natural, we will use the same symbol $\pi_{\leq 0}(\Def({\mathfrak g}))$ to denote both the set $\pi_{\leq 0}(\Def({\mathfrak g}))$ and the pointed set $\pi_{\leq 0}(\Def({\mathfrak g});0)$.
\vskip .5 cm  

It is important to remark that the functors of the form $\Def(\mathfrak{g})$ are very special ones among all formal $\infty$-groupoids. To begin with, $\Def(\mathfrak{g})(\mathbb{K})=\{0\}$ and so, in particular, $\Def(\mathfrak{g})(\mathbb{K})$ is a homotopically trivial $\infty$-groupoid.
Another characterzing property of the functors of the form $\Def(\mathfrak{g})$  among formal $\infty$-groupoids is that, under suitable assumptions, they commute with homotopy pullbacks; see \cite{lurie-moduli, pridham} for a precise statement. In other words, if we call ``formal moduli problems'' those formal $\infty$-groupoids which satisfy the two conditions we have just observed for $\Def(\mathfrak{g})$, what we are saying is that $\Def$ is actually a functor
\[ \Def: {\bf DGLA} \to {\text{\bf Formal moduli problems}}.\]
And a very good reason for working with $\infty$-groupoids valued deformation functors rather than with 
their apparently handier set-valued or groupoid-valued versions is the following remarkable result, which allows one to move homotopy constructions back and forth between dglas and formal moduli problems.
\begin{theorem*}[Pridham-Lurie]
The functor $\Def: {\bf DGLA} \to {\text{\bf Formal moduli problems}}$ is an equivalence of 
$(\infty,1)$-categories.
\end{theorem*}
Here the $(\infty,1)$-category structures involved are the most natural ones, and they are both induced by standard model category structures. Namely, on the category of dglas one takes surjective morphisms as fibrations and quasi-isomorphisms as weak equivalences, just as in the case of differential complexes, whereas the model category structure on the right hand side is induced by the standard model category structure on Kan complexes as a subcategory of simplicial sets.
A proof of the above equivalence can be found in \cite{lurie-moduli, pridham}.
\par
\vskip .5 cm
We will often identify a dgla $\mathfrak{g}$ with the functor ${\bf dgArt}\to {\bf nilpotent\,\, DGLA}$ it defines by the rule $A\mapsto \mathfrak{g}\otimes \mathfrak{m}_A$. With this in mind, we will occasionally apply constructions that generally only
make sense for nilpotent dglas (such as $\exp$) to arbitrary dglas. What we mean in these cases is that  the construction is applied not
to a single dgla, but to the functor from ${\bf dgArt}$ to nilpotent dglas it defines. The same kind of consideration applies to our somehow colloquial use of the expression ``$\infty$-groupoid'' in the following sections; namely, by that we will occasionally mean ``formal $\infty$-groupoid'', or even ``formal stack in $\infty$-groupoids''. The precise meaning to be given to ``$\infty$-groupoid''  will always be clear from the context. 

\section{Tangent spaces and obstructions}
If $\mathbf{X}$ is a formal moduli problem, then the simplicial set $\mathbf{X}(\mathbb{K}[\epsilon]/(\epsilon^2))$ has a natural structure of simplicial vector space, and so, via the Dold-Kan correspondence, it is equivalent to the datum of a chain complex: the tangent complex $T\mathbf{X}$ of $\mathbf{X}$. Passing from $\mathbf{X}$ to the associated classical moduli problem $\pi_{\leq 0}\mathbf{X}$, the only datum we read of the tangent complex is its homotopy class, i.e., since we are working on a field, its cohomology. In particular, we have a natural isomorphism
\[
T\pi_{\leq 0}\Def \xrightarrow{\sim} H^1
\]
of functors $\mathbf{DGLA}\to \text{\bf Vector spaces}$ between the tangent space to the classical moduli problem associated to a dgla and the first cohomology group of the dgla seen as a cochain complex. Let us rephrase this in a more explicit form. As we noticed in the previous section, $\pi_{\leq 0}\Def({\mathfrak g})$ is the functor of Artin rings $
A\mapsto \MC({\mathfrak g}\otimes {\mathfrak m}_A)\bigl/{\rm gauge}$, hence 
\[
T\pi_{\leq 0}\Def (\mathfrak{g})=\MC\left({\mathfrak g}\otimes (\epsilon)/(\epsilon^2)\right)\bigl/{\rm gauge}\cong H^1(\mathfrak{g}).
\]
This isomorphism is natural. Namely, given a morphism $\varphi\colon \mathfrak{g}\to \mathfrak{h}$ of dglas, let us write $\Phi$ for the induced morphism of classical moduli problems,
\[
\Phi=\pi_{\leq 0}\Def(\varphi)\colon \pi_{\leq 0}\Def(\mathfrak{g})\to \pi_{\leq 0}\Def(\mathfrak{h}).
\]
Then the differential of $\Phi$,
\[
d\Phi\colon T\pi_{\leq 0}\Def(\mathfrak{g})\to T\pi_{\leq 0}\Def(\mathfrak{h})
\]
is naturally identified with 
\[
H^1(\varphi):H^1(\mathfrak{g})\to H^1(\mathfrak{h}).
\]
The second cohomology group $H^2$ defines a natural 
obstruction theory for $\pi_{\leq 0}\Def$, i.e., obstructions for the classical moduli problem $\pi_{\leq 0}\Def(\mathfrak{g})$ are naturally identified with elements in $H^2(\mathfrak{g})$, see \cite{manetti-deformations}. Note that this does not mean that each element in $H^2(\mathfrak{g})$ represents an obstruction: one can have dglas with nontrivial $H^2$ governing unobstructed deformation problems. The naturally of the obstruction theory given by the second cohomology groups means that, if $\varphi\colon \mathfrak{g}\to \mathfrak{h}$ is a morphism of of dglas, the induced morphism in cohomology,
 \[
H^2(\varphi):H^2(\mathfrak{g})\to H^2(\mathfrak{h}),
\]
maps obstructions for the classical moduli problem $\pi_{\leq 0}\Def(\mathfrak{g})$ to obstructions for the classical moduli problem $\pi_{\leq 0}\Def(\mathfrak{h})$. In particular, if the moduli problem $\pi_{\leq 0}\Def(\mathfrak{h})$ is unobstructed (e.g., if the functor $\pi_{\leq 0}\Def(\mathfrak{h})$ is smooth), then 
\[
\mathrm{Obstructions}\left(\pi_{\leq 0}\Def(\mathfrak{g})\right)\subseteq \ker\left(H^2(\varphi)\colon H^2(\mathfrak{g})\to H^2(\mathfrak{h})\right). 
\]

\section{Homotopy vs. gauge equivalent morphisms of dglas (with a detour into $L_\infty$-morphisms)}

 Let ${\mathfrak g}$ and ${\mathfrak h}$ be two dglas. 
The hom-space $\Hom_\infty(\mathfrak{g},\mathfrak{h})$ of morphisms
between $\mathfrak{g}$ and $\mathfrak{h}$ in the $(\infty,1)$-category of dglas is conveniently modelled as the simplicial set $\MC(\underline{\Hom}({\mathfrak g},{\mathfrak h})\otimes\Omega_\bullet)$, where $\underline{\Hom}({\mathfrak g},{\mathfrak h})$ is the
Chevalley-Eilenberg-type dgla associated with the pair $(\mathfrak{g},\mathfrak{h})$. It is given as the total dgla of the bigraded dgla
\[
\underline{\Hom}^{p,q}({\mathfrak g},{\mathfrak h})=\Hom_{{\mathbb Z}-{\bf Vect}}(\wedge^q{\mathfrak g},{\mathfrak h}[p])=\Hom^p(\wedge^q{\mathfrak g},{\mathfrak h}),
\]
endowed with the Lie bracket
\[
[\,,\,]_{\underline{\Hom}}\colon \underline{\Hom}^{p_1,q_1}({\mathfrak g},{\mathfrak h})\otimes \underline{\Hom}^{p_2,q_2}({\mathfrak g},{\mathfrak h})\to \underline{\Hom}^{p_1+p_2,q_1+q_2}({\mathfrak g},{\mathfrak h})
\]
defined by
\begin{align*}
[f,g]^{}_{\underline{\Hom}}&(\gamma_1^{}\wedge\cdots\wedge\gamma_{q_1+q_2}^{})=\\
&\hskip-1em=\sum_{\sigma\in{\rm Sh}(q_1,q_2)}\pm[f(
\gamma_{\sigma(1)}\wedge\cdots \wedge\gamma_{\sigma(q_1)}),
g(
\gamma_{\sigma(q_1+1)}\wedge\cdots \wedge\gamma_{\sigma(q_1+q_2)})]^{}_{\mathfrak h},
\end{align*}
with $\sigma$ ranging in the set of $(q_1,q_2)$-unshuffles and and $\pm$ standing for the Koszul sign, 
and with the differentials
\[
d_{1,0}^{}\colon
\underline{\Hom}^{p,q}({\mathfrak g},{\mathfrak h})\to \underline{\Hom}^{p+1,q}({\mathfrak g},{\mathfrak h})\]
and
\[
d_{0,1}^{}\colon
\underline{\Hom}^{p,q}({\mathfrak g},{\mathfrak h})\to \underline{\Hom}^{p,q+1}({\mathfrak g},{\mathfrak h})\]
given by
\[
(d_{1,0}^{}{f})(\gamma_1^{}\wedge\cdots\wedge\gamma_q^{})=d_{\mathfrak h}(f(\gamma_1^{}\wedge\cdots\wedge\gamma_q^{}))+\sum_i
\pm f(\gamma_1\wedge\cdots
\wedge d_{\mathfrak g}\gamma_i\wedge
\cdots\wedge\gamma_{q+1}^{})
\]
and 
\[
(d_{0,1}^{}f)(\gamma_1\wedge\cdots\wedge\gamma_{q+1})=
\sum_{i<j}\pm f([\gamma_i,\gamma_j]^{}_{\mathfrak g}\wedge
\gamma_1\wedge\cdots
\wedge\widehat{\gamma_i}\wedge\cdots\wedge\widehat{\gamma_j}\wedge
\cdots\wedge\gamma_{q+1}^{}).
\]
An explicit determination for the signs in the above formulas can be found, e.g, in \cite{lada-markl,schuhmacher}.
These operations are best seen pictorially:
\[
\left[\quad
\begin{xy}
,(-5,-4.5);(0,0)*{\circ}**\dir{-}
,(-2,-6);(0,0)*{\circ}**\dir{-}
,(2,-6);(0,0)*{\circ}**\dir{-}
,(5,-4.5);(0,0)*{\circ}**\dir{-}
,(0,0)*{\circ};(0,6)**\dir{-}?>*\dir{>}
,(-2.4,0.5)*{f}
\end{xy}\,,\,\,
\begin{xy}
,(-3,-5.5);(0,0)*{\circ}**\dir{-}
,(1,-6);(0,0)*{\circ}**\dir{-}
,(4,-6);(0,0)*{\circ}**\dir{-}
,(2.4,0.5)*{g}
,(0,0)*{\circ};(0,6)**\dir{-}?>*\dir{>}
\end{xy}\quad
\right]_{\underline{\Hom}}\phantom{i}=\phantom{m}
\begin{xy}
,(-5,-6.5);(0,-2)*{\circ}**\dir{-}
,(-2,-8);(0,-2)*{\circ}**\dir{-}
,(2,-8);(0,-2)*{\circ}**\dir{-}
,(5,-6.5);(0,-2)*{\circ}**\dir{-}
,(0,-2)*{\circ};(6,2.5)*{\bullet}**\dir{-}
,(-2.4,-1.5)*{f}
,(8,-7.5);(11,-2)*{\circ}**\dir{-}
,(12,-8);(11,-2)*{\circ}**\dir{-}
,(16,-8);(11,-2)*{\circ}**\dir{-}
,(13.6,-1.5)*{g}
,(11,-2)*{\circ};(6,2.5)*{\bullet}**\dir{-}
,(6,2.5)*{\bullet};(6,8)**\dir{-}?>*\dir{>}
,(9.5,3.5)*{\scriptstyle{[\,,\,]_{\mathfrak{h}}^{}}}
\end{xy}\qquad;
\]
\vskip .5 cm
\[
d_{1,0}\left(\,
\begin{xy}
,(-5,-4.5);(0,0)*{\circ}**\dir{-}
,(-2,-6);(0,0)*{\circ}**\dir{-}
,(2,-6);(0,0)*{\circ}**\dir{-}
,(5,-4.5);(0,0)*{\circ}**\dir{-}
,(0,0)*{\circ};(0,6)**\dir{-}?>*\dir{>}
,(-2.4,0.5)*{f}
\end{xy}\,
\right)\phantom{i}=\phantom{i}
\begin{xy}
,(-7,-6.5);(0,0)*{\circ}**\dir{-}
,(-2,-7);(0,0)*{\circ}**\dir{-}
,(2,-7);(0,0)*{\circ}**\dir{-}
,(7,-6.5);(0,0)*{\circ}**\dir{-}
,(0,0)*{\circ};(0,7)**\dir{-}?>*\dir{>}
,(-2.4,0.5)*{f}
,(0,3)*{\bullet}
,(2.8,3.5)*{\scriptstyle{d_{\mathfrak{h}}}}
\end{xy}
\phantom{i}+\phantom{i}
\begin{xy}
,(-7,-6.5);(0,0)*{\circ}**\dir{-}
,(-2,-7);(0,0)*{\circ}**\dir{-}
,(2,-7);(0,0)*{\circ}**\dir{-}
,(7,-6.5);(0,0)*{\circ}**\dir{-}
,(0,0)*{\circ};(0,6)**\dir{-}?>*\dir{>}
,(-2.4,0.5)*{f}
,(-4,-3.7)*{\bullet}
,(-6.3,-2.4)*{\scriptstyle{d_{\mathfrak{g}}}}
\end{xy}
\qquad;\qquad
d_{0,1}
\left(\,
\begin{xy}
,(-5,-4.5);(0,0)*{\circ}**\dir{-}
,(-2,-6);(0,0)*{\circ}**\dir{-}
,(2,-6);(0,0)*{\circ}**\dir{-}
,(5,-4.5);(0,0)*{\circ}**\dir{-}
,(0,0)*{\circ};(0,6)**\dir{-}?>*\dir{>}
,(-2.4,0.5)*{f}
\end{xy}\,
\right)\phantom{i}=\phantom{i}
\begin{xy}
,(-4,-3.7);(0,0)*{\circ}**\dir{-}
,(-4,-3.7);(-5.8,-7)**\dir{-}
,(-4,-3.7);(-2.8,-7)**\dir{-}
,(-1,-7);(0,0)*{\circ}**\dir{-}
,(3,-7);(0,0)*{\circ}**\dir{-}
,(7,-6.5);(0,0)*{\circ}**\dir{-}
,(0,0)*{\circ};(0,6)**\dir{-}?>*\dir{>}
,(-2.4,0.5)*{f}
,(-4,-3.7)*{\bullet}
,(-7.5,-2.5)*{\scriptstyle{[\,,\,]_{\mathfrak{g}}}}
\end{xy}\quad.\]
\vskip.2 cm
Note that $0$-simplices in $\MC(\underline{\Hom}({\mathfrak g},{\mathfrak h})\otimes\Omega_\bullet)$ are those degree 1 elements $f$ in $\underline{\Hom}({\mathfrak g},{\mathfrak h})$ such that
\begin{align*}
&d_{\mathfrak h}f_n(\gamma_1^{}\wedge\cdots\wedge\gamma_n^{})+\frac{1}{2}\!\!\!\sum_{\substack{q_1+q_2=n\\ \sigma\in{\rm Sh}(q_1,q_2)}}\!\!\!\!\!\!\pm[f_{q_1}(
\gamma_{\sigma(1)}\wedge\cdots \wedge\gamma_{\sigma(q_1)}),
f_{q_2}(
\gamma_{\sigma(q_1+1)}\wedge\cdots \wedge\gamma_{\sigma(q_1+q_2)})]^{}_{\mathfrak h}\\
&=\sum_i
\pm f_n(\gamma_1\wedge\cdots
\wedge d_{\mathfrak g}\gamma_i\wedge
\cdots\wedge\gamma_{n}^{})\\
&\qquad\qquad+
\sum_{i<j}\pm f_{n-1}([\gamma_i,\gamma_j]^{}_{\mathfrak g}\wedge
\gamma_1\wedge\cdots
\wedge\widehat{\gamma_i}\wedge\cdots\wedge\widehat{\gamma_j}\wedge
\cdots\wedge\gamma_{q+1}^{}),
\end{align*}
where we have written $f=\sum_{n\geq 1}f_n$ with $f_n$ in $\underline{\Hom}^{1-n,n}({\mathfrak g},{\mathfrak h})$.
Here one recognizes the explicit equation characterizing $L_\infty$-morphisms between $\mathfrak{g}$ and $\mathfrak{h}$ one finds  in, e.g., \cite{lada-stasheff,kontsevich}.

%
 Next, as we remarked in the previous section, the set $\pi_{\leq 0}(\MC(\underline{\Hom}({\mathfrak g},{\mathfrak h})\otimes\Omega_\bullet))$ is isomorphic to the quotient $\MC(\underline{\Hom}({\mathfrak g},{\mathfrak h}))/{\rm gauge}$. Thus we obtain the following proposition.
\begin{proposition*}
Let $f,g:{\mathfrak g} \to {\mathfrak h}$ be two $L_\infty$-morphisms of dglas. Then $f$ and $g$ are gauge equivalent in $\MC(\underline{\Hom}({\mathfrak g},{\mathfrak h}))$ if and only if $f$ and $g$ represent the same morphism in the homotopy category of dglas.
\end{proposition*}

\vskip .3 cm
It should be remarked that the above construction is an instance of a more general phenomenon: 
if $\mathcal{O}$ is a differential graded operad, $A$ is an $\mathcal{O}$-algebra, and $B$ is a (differential graded) cocommutative coalgebra, then the space of linear mappings from $B$ to $A$ has a natural $\mathcal{O}$-algebra structure, see \cite{dolgushev}.
%
%
 Also, the analogue of the above proposition holds, more in general, for the homotopy category of ${\mathcal O}$-algebras, see \cite{merkulov-vallette, pridham2}. We thank Jonathan Pridham and Bruno Vallette for this remark.

\section{Cartan homotopies appear}
Let now ${\mathfrak g}$ and ${\mathfrak h}$ be dglas and $\boldsymbol{i}\colon {\mathfrak g}\to {\mathfrak h}[-1]$ be a morphism of graded vector spaces. Let $\lie=d_{1,0}\boldsymbol{i}$; that is, $\lie_a= d_{\mathfrak h} \boldsymbol{i}_a + \boldsymbol{i}_{d_{\mathfrak g}a}$ for any $a\in{\mathfrak g}$. The morphism $\boldsymbol{i}$ is called a Cartan homotopy\footnote{Up to our knowledge, this terminology has been introduced in \cite{fiorenza-manetti, fiorenza-manetti2}.} if it satisfies the two conditions
\[   \boldsymbol{i}_{[a,b]_{\mathfrak g}}= [\boldsymbol{i}_a, \lie_b]_{\mathfrak h}\qquad \text{and} \qquad [\boldsymbol{i}_a, \boldsymbol{i}_b]_{\mathfrak h}=0,\qquad \mbox{for all} \ a, b \in {\mathfrak g}. \]
This name has an evident geometric origin: if ${\mathcal T}_X$ is the tangent sheaf of a smooth manifold $X$ and $\Omega^*_{X}$ is the sheaf of complexes of differential forms, then the contraction of differential forms with vector fields is a Cartan homotopy
\[
\boldsymbol{i}\colon {\mathcal T}_X\to {\mathcal E}nd^*(\Omega^*_X)[-1].
\]
In this case, $\lie_a$ is the Lie derivative along the vector field $a$, and the conditions $\boldsymbol{i}_{[a,b]}= [\boldsymbol{i}_a, \lie_b]$ and $[\boldsymbol{i}_a, \boldsymbol{i}_b]=0$, together with the defining equation $\lie_a=[d_{\Omega^*_X},\boldsymbol{i}_a]$ and with the equations $\lie_{[a,b]}=[\lie_a,\lie_b]$ and $[d_{\Omega^*_X},\lie_a]=0$ expressing the fact that $\lie\colon {\mathcal T}_X\to {\mathcal E}nd^*(\Omega^*_X)$ is a dgla morphism, are nothing but the well-known Cartan identities involving contractions and Lie derivatives.
\par
It is a straightforward computation to see that, if $\boldsymbol{i}$ is a Cartan homotopy, then the degree zero morphism of graded vector spaces $\lie\colon \mathfrak{g}\to \mathfrak{h}$ is actually a dgla morphism. Actually much more is true: the dgla morphism $\lie$ is gauge equivalent to zero, and $e^{\boldsymbol{i}}$ is a homotopy between $\lie$ and $0$. To see this, just compute 
%
\begin{align*}e^{-\boldsymbol{i}}*0 &= \sum_{n=0}^{+\infty} \frac{{({\rm ad}_{-\boldsymbol{i}})}^n}{(n+1)!}\ (d_{\underline{\Hom}}\boldsymbol{i})= \sum_{n=0}^{+\infty} \frac{{({\rm ad}_{-\boldsymbol{i}})}^n}{(n+1)!}\ (\lie+\boldsymbol{i}_{[\,,\,]_{\mathfrak g}}).\end{align*}
Then one sees that  the $(0,1)$-component of $e^{-\boldsymbol{i}}*0$ 
is just $\lie$; the $(-1,2)$-component is
\[
\boldsymbol{i}_{[\,,\,]_{\mathfrak g}}-\frac{1}{2}[\boldsymbol{i},\lie]_{\underline{\Hom}}
\] 
and, for $n\geq 3$, the $(1-n,n)$-component has two contributions, one of the form \allowbreak $[\boldsymbol{i},[\boldsymbol{i},\cdots,[\boldsymbol{i}, \lie]_{\underline{\Hom}}\cdots]_{\underline{\Hom}}]_{\underline{\Hom}}$ and the other of the form $[\boldsymbol{i},[\boldsymbol{i},\cdots,[\boldsymbol{i}, \boldsymbol{i}_{[\,,\,]_{\mathfrak g}}]_{\underline{\Hom}}\cdots]_{\underline{\Hom}}]_{\underline{\Hom}}$. Threfore, if $\boldsymbol{i}$ is a Cartan homotopy, then all the nonlinear components of $e^{-\boldsymbol{i}}*0$ vanish and one finds $e^{-\boldsymbol{i}}*0=\lie$, i.e.,
\[
e^{\boldsymbol{i}}*\lie=0.
\]
The above discussion can be summarized as follows.
\begin{proposition*} 
Let ${\mathfrak g}$ and ${\mathfrak h}$ be two dglas. If $\boldsymbol{i}\colon {\mathfrak g}\to{\mathfrak h}[-1]$ is a Cartan homotopy, then $\lie=d_{1,0}\boldsymbol{i}\colon {\mathfrak g}\to{\mathfrak h}$ is a dgla morphism gauge equivalent to the zero morphism via the gauge action of $e^{\boldsymbol{i}}$.
\end{proposition*}

\section{Homotopy fibers (and the associated exact sequence)}
Let now $\boldsymbol{i}\colon {\mathfrak g}\to{\mathfrak h}[-1]$ be a Cartan homotopy and $\lie\colon {\mathfrak g}\to{\mathfrak h}$ be the associated dgla morphism. Then, the equation $e^{\boldsymbol{i}}*\boldsymbol{l}=0$ implies that, for any subdgla ${\mathfrak n}$ of ${\mathfrak h}$ containing the image of $\lie$, the morphism $\lie\colon {\mathfrak g}\to{\mathfrak n}$ equalizes the diagram 
$\xymatrix{ {\mathfrak n} \ar@<2pt>[r]^{{\rm incl.}}\ar@<-2pt>[r]_{0} & {\mathfrak h}}$
 up to a homotopy provided by the gauge action of $e^{\boldsymbol{i}}$. Hence we have a morphism to the homotopy limit:
\[ {\mathfrak g} \xrightarrow{(\lie, e^{\boldsymbol{i}})}\holim \left(\xymatrix{ {\mathfrak n} \ar@<2pt>[r]^{{\rm incl.}}\ar@<-2pt>[r]_{0} & {\mathfrak h}}\right). \]
Taking Def's we obtain a natural transformation of $\infty$-groupoid valued functors:
\[ \Def({\mathfrak g})\xrightarrow{(\lie, e^{\boldsymbol{i}})}\holim \left(\xymatrix{ \Def({\mathfrak n}) \ar@<2pt>[r]^{\Def_{\rm incl.}}\ar@<-2pt>[r]_{\Def_0} &\Def({\mathfrak h})}\right). 
\]
The map $\Def_0\colon \Def({\mathfrak n})\to \Def({\mathfrak h})$ is the constant map to the distinguished point $0$ in $\Def({\mathfrak h})$; therefore, the homotopy limit above is the homotopy fiber of 
$\Def_{\rm incl.}\colon \Def({\mathfrak n})\to \Def({\mathfrak h})$ over the point $0$, and we obtain a natural transformation
\[ \Def({\mathfrak g})\xrightarrow{(\lie, e^{\boldsymbol{i}})} {\rm hoDef}^{-1}_{\rm incl.}(0), 
\]
which at the zeroth level gives a natural transformation of {\bf Set}-valued deformation functors
\[
{\mathcal P}\colon \pi_{\leq 0}\Def({\mathfrak g})\to\pi_{\leq 0}{\rm hoDef}^{-1}_{\rm incl.}(0).
\]
The differential of ${\mathcal P}$ is easily computed: it is the linear map 
\[ 
H^1({\mathfrak g}) \xrightarrow{H^1((\lie, e^{\boldsymbol{i}}))}
H^1(\holim \left(\xymatrix{ {\mathfrak n} \ar@<2pt>[r]^{{\rm incl.}}\ar@<-2pt>[r]_{0} & {\mathfrak h}}\right)).
 \]
 Since the model category structure on dglas is the same as on differential complexes, we can compute the $H^1$ on the right hand side by taking the holimit in complexes. Then the natural quasi-isomorphism $\holim (\xymatrix{ {\mathfrak n} \ar@<2pt>[r]^{{\rm incl.}}\ar@<-2pt>[r]_{0} & {\mathfrak h}})\simeq ({\mathfrak h}/{\mathfrak n})[-1]$ tells us that the differential of ${\mathcal P}$ is just the  map 
 \[
 H^1({\boldsymbol{i}})\colon H^1({\mathfrak g}) \to H^0({\mathfrak h}/{\mathfrak n})
 \]
induced by the morphism of complexes ${\boldsymbol{i}}\colon {\mathfrak g} \to ({\mathfrak h}/{\mathfrak n})[-1]$. Also, the map
\[
 H^2({\boldsymbol{i}})\colon H^2({\mathfrak g}) \to H^1({\mathfrak h}/{\mathfrak n})
 \]
 maps the obstruction space of $\pi_{\leq 0}\Def({\mathfrak g})$ (as a subspace of  $H^2({\mathfrak g})$) to the obstruction space of  $\pi_{\leq 0}{\rm hoDef}^{-1}_{\rm incl.}(0)$ (as a subspace of $H^1({\mathfrak h}/{\mathfrak n})$). In particular, if $\pi_{\leq 0}{\rm hoDef}^{-1}_{\rm incl.}(0)$ is smooth, and therefore unobstructed, the obstructions of the deformation functor $\pi_{\leq 0}\Def({\mathfrak g})$ are contained in the kernel of the map $ H^2({\boldsymbol{i}})\colon H^2({\mathfrak g}) \to H^1({\mathfrak h}/{\mathfrak n})$.
 \vskip .8 cm
 
To investigate the geometrical aspects of the map $\mathcal{P}$, note that, by
looking at $\pi_{\leq 0}{\rm hoDef}^{-1}_{\rm incl.}(0)$ as a pointed set, it nicely fits into the homotopy exact sequence
\[
\pi_1(\Def({\mathfrak n}); 0) \xrightarrow{\Def_{{\rm incl.}*}} \pi_1(\Def({\mathfrak h}); 0) \to \pi_{0}({\rm hoDef}^{-1}_{\rm incl.}(0);0) \to  \pi_0(\Def({\mathfrak n}); 0), 
\]
so we get a canonical isomorphism between the preimage of the distinguished point $0$ under the map $\pi_{0}({\rm hoDef}^{-1}_{\rm incl.}(0);0) \to  \pi_0(\Def({\mathfrak n}); 0)$ and the quotient set
\[
\frac{ \pi_1(\Def({\mathfrak h});0)}{{\Def_{{\rm incl.}*}}\pi_1(\Def({\mathfrak n});0)}.
\]
The group $\pi_1(\Def({\mathfrak h});0)$ is the group of automorphisms of $0$ in the groupoid $\pi_{\leq 1}(\Def({\mathfrak h}))$. We have already remarked that this groupoid is not equivalent to the Deligne groupoid of ${\mathfrak h}$, i.e., the action groupoid for the gauge action of $\exp({\mathfrak h}^0\otimes {\mathfrak m}_A)$ on $\MC({\mathfrak h}\otimes {\mathfrak m}_A)$, since the irrelevant stabilizer 
\[ {\rm Stab}(x)= \{dh+[x,h] \mid  h \in \mathfrak{h}^{-1}\otimes
 {\mathfrak m}_A\}\subseteq \{ a \in \mathfrak{h}^0\otimes \mathfrak{m}_A \mid e^a*x=x \}
 \]
 of a Maurer-Cartan element $x$ may be nontrivial. However, the group  $\pi_1(\Def({\mathfrak h});0)$ only sees the connected component of $0$, and on this connected component the irrelevant stabilizers are trivial as soon as the differential of the dgla ${\mathfrak h}$ vanishes on $\mathfrak{h}^{-1}$. This immediately follows from noticing that irrelevant stabilizers of gauge equivalent Maurer-Cartan elements are conjugate subgroups of $\exp(\mathfrak{h}^0\otimes\mathfrak{m}_A)$, see, e.g., \cite{manetti-obstructions}. In particular, if ${\mathfrak h}$ is a graded Lie algebra (which we can consider as a dgla with trivial differential), then 
$\pi_1(\Def({\mathfrak h});0)\simeq \exp({\mathfrak h}^0)$, 
where ${\mathfrak h}^0$ denotes the degree zero component of ${\mathfrak h}$. Similarly, since ${\mathfrak n}$ is a subdgla of ${\mathfrak h}$, one has $\pi_1(\Def({\mathfrak n});0)\simeq \exp({\mathfrak n}^0)$, and the group homomorphism $\Def_{{\rm incl.}*}$ is just the inclusion. Moreover, when ${\mathfrak h}$ has trivial differential, $e^\alpha*0=0$ for any $\alpha$ in $\mathfrak{h}^0\otimes\mathfrak{m}_A$. Therefore, the composition $\pi_{\leq 0}\Def({\mathfrak g})\to \pi_{\leq 0}{\rm hoDef}^{-1}_{\rm incl.}(0)\to \pi_{\leq 0}\Def({\mathfrak n})$
maps the whole of  $\pi_{\leq 0}\Def({\mathfrak g})$ onto the distinguished point $0$ of  $\pi_{\leq 0}\Def({\mathfrak n})$ and so it induces a natural map $\pi_{\leq 0}\Def({\mathfrak g})\to \exp({\mathfrak h}^0)/\exp({\mathfrak n}^0)$. This is nothing but the natural map
\[
e^{\boldsymbol{i}}\colon \pi_{\leq 0}\Def({\mathfrak g})\to \exp({\mathfrak h}^0)/\exp({\mathfrak n}^0)
\]
which sends a Maurer-Cartan element $\xi\in {\mathfrak g}^1\otimes {\mathfrak m}_A$ to $e^{\boldsymbol{i}_\xi}\mod \exp({\mathfrak n}^0\otimes {\mathfrak m}_A)$. A particularly interesting case is when the pair $({\mathfrak h},{\mathfrak n})$ is formal,\footnote{We are not sure whether this terminology is a standard one} i.e., if the inclusion of ${\mathfrak n}$ in ${\mathfrak h}$ induces an inclusion in cohomology and the two inclusions $H^*({\mathfrak n})\hookrightarrow  H^*({\mathfrak h})$ and ${\mathfrak n}\hookrightarrow {\mathfrak h}$ are homotopy equivalent.
 Indeed, in this case the pair $(\Def({\mathfrak h}),\Def({\mathfrak n}))$ will be equivalent to the pair  $(\Def(H^*({\mathfrak h})),\Def(H^*({\mathfrak h})))$ and there will be an induced isomorphism between $\pi_1(\Def({\mathfrak h});0)/{\Def_{{\rm incl.}*}}\pi_1(\Def({\mathfrak n});0)$ and the smooth homogeneous space $\exp(H^0({\mathfrak h}))/\exp(H^0({\mathfrak n}))$. We can summarize the results described in this section as follows:
 \begin{proposition*}
 
Let $\boldsymbol{i}\colon {\mathfrak g}\to{\mathfrak h}[-1]$ be a Cartan homotopy, let $\lie\colon {\mathfrak g}\to{\mathfrak h}$ be the associated dgla morphism, and let ${\mathfrak n}$  be a subdgla of ${\mathfrak h}$ containing the image of $\lie$. Then, if the pair $({\mathfrak h},{\mathfrak n})$ is formal, we have a natural transformation\footnote{This natural transformation is not canonical: it depends on the choice of a quasi isomorphism $({\mathfrak h},{\mathfrak n})\simeq(H^*({\mathfrak h}),H^*({\mathfrak n}))$. Also note that the tangent space at $0$ on the right hand side is $H^0({\mathfrak h})/H^0({\mathfrak n})$; this is only apparently in contrast with the general result mentioned above that the tangent space at $0$ to $\pi_{\leq 0}{\rm hoDef}^{-1}_{\rm incl.}(0)$ is $H^0({\mathfrak h}/{\mathfrak n})$. Indeed, when $({\mathfrak h},{\mathfrak n})$ is a formal pair, the two vector spaces  $H^0({\mathfrak h})/H^0({\mathfrak n})$ and $H^0({\mathfrak h}/{\mathfrak n})$ are (non canonically) isomorphic.} of {\bf Set}-valued deformation functors
\[
{\mathcal P}\colon \pi_{\leq 0}(\Def({\mathfrak g})) \to \exp(H^0({\mathfrak h}))/\exp(H^0({\mathfrak n}))
\]
induced by the dgla map
\[ {\mathfrak g} \xrightarrow{(\lie, e^{\boldsymbol{i}})}\holim \left(\xymatrix{ {\mathfrak n} \ar@<2pt>[r]^{{\rm incl.}}\ar@<-2pt>[r]_{0} & {\mathfrak h}}\right). \]
In particular, since $\exp(H^0({\mathfrak h}))/\exp(H^0({\mathfrak n}))$ is smooth, the obstructions of the {\bf Set}-valued deformation functor 
 $\pi_{\leq 0}(\Def({\mathfrak g}); 0)$ are contained in the kernel of the map  $H^2({\boldsymbol{i}})\colon H^2({\mathfrak g}) \to H^1({\mathfrak h}/{\mathfrak n})$.
 \end{proposition*}
 This result can be nicely refined, by showing how the main result from \cite{iacono-manetti} naturally fits into the discussion above. We have:
 \begin{proposition*}
 Let $({\mathfrak h},{\mathfrak n})$ be a formal pair of dglas. Then, the dgla $\holim \left(\xymatrix{ {\mathfrak n} \ar@<2pt>[r]^{{\rm incl.}}\ar@<-2pt>[r]_{0} & {\mathfrak h}}\right)$ is quasi-abelian. In particular there is a (non-canonical) quasi-isomorphism of dglas between $\holim \left(\xymatrix{ {\mathfrak n} \ar@<2pt>[r]^{{\rm incl.}}\ar@<-2pt>[r]_{0} & {\mathfrak h}}\right)$ and the abelian dgla obtained by endowing the complex $({\mathfrak h}/{\mathfrak n})[-1]$ with the trivial bracket.
 \end{proposition*}
 To see this, notice that, since by hypothesis the inclusion ${\mathfrak n}\hookrightarrow {\mathfrak h}$ induces an inclusion $H^*({\mathfrak n})\hookrightarrow H^*({\mathfrak h})$, the projection ${\mathfrak h}[-1]\to {\mathfrak h}/{\mathfrak n}[-1]$ admits a section $\boldsymbol{i}$ which is a morphism of complexes. Denote by ${\mathfrak g}$ the dgla obtained from the complex ${\mathfrak h}/{\mathfrak n}[-1]$ by endowing it with the trivial bracket. Then, the map of graded vector spaces $\boldsymbol{i}\colon {\mathfrak g}\to {\mathfrak h}[-1]$ is a Cartan homotopy whose associated dgla morphism is the zero map $0\colon {\mathfrak g}\to{\mathfrak h}$. Therefore we have a dgla map
 \[ ({\mathfrak h}/{\mathfrak n})[-1] \xrightarrow{(0, e^{\boldsymbol{i}})}\holim \left(\xymatrix{ {\mathfrak n} \ar@<2pt>[r]^{{\rm incl.}}\ar@<-2pt>[r]_{0} & {\mathfrak h}}\right). \]
 Since $\boldsymbol{i}$ is a section to ${\mathfrak h}[-1]\to {\mathfrak h}/{\mathfrak n}[-1]$,  the map in cohomology
 \[
 H^*({\mathfrak h}/{\mathfrak n})[-1] \xrightarrow{H^*(0,e^{\boldsymbol{i}})}H^*(\holim \left(\xymatrix{ {\mathfrak n} \ar@<2pt>[r]^{{\rm incl.}}\ar@<-2pt>[r]_{0} & {\mathfrak h}}\right)) \]
 is identified with the identity of $H^*({\mathfrak h}/{\mathfrak n})[-1]$ by the the natural quasi-isomorphism of complexes $\holim (\xymatrix{ {\mathfrak n} \ar@<2pt>[r]^{{\rm incl.}}\ar@<-2pt>[r]_{0} & {\mathfrak h}})\xrightarrow{\sim} ({\mathfrak h}/{\mathfrak n})[-1]$.

\section{From local to global, and classical (and generalized) periods}
Assume now $\mathbb{K}$ is algebraically closed. Let $X$ be a smooth projective manifold, and let ${\mathcal T}_X$ and $\Omega_X^*$ be the tangent sheaf and the sheaf of differential forms on $X$, respectively. The sheaf of complexes $(\Omega_X^*,d_{\Omega_X^*})$ is naturally filtered by setting $F^p\Omega_X^*=\oplus_{i\geq p}\Omega^i_X$. Finally, let ${\mathcal E}nd^*(\Omega_X^*)$ be the endomorphism sheaf of $\Omega_X^*$ and ${\mathcal E}nd^{\geq 0}(\Omega_X^*)$ be the subsheaf consisting of nonnegative degree elements. Note that ${\mathcal E}nd^{\geq 0}(\Omega_X^*)$ is a subdgla of ${\mathcal E}nd^{*}(\Omega_X^*)$, and can be seen as the subdgla of endomorphisms preserving the filtration on $\Omega_X^*$.

Recall that the prototypical example of Cartan homotopy was the contraction of differential forms with vector fields ${\boldsymbol{i}}: \T_X \to \Eps nd^*(\Omega_X^*)[-1]$; the corresponding 
dgla morphism  is $a\mapsto \lie_a$, where $\lie_a$ the Lie derivative along $a$. Explicitly, $\lie_a= d_{\Omega_X^*}\circ {\boldsymbol{i}}_a+{\boldsymbol{i}}_{a}\circ d_{\Omega_X^*}$, and so $\lie_a$ preserves the filtration. Therefore, we have a natural transformation\footnote{Of what? The correct answer would be of $\infty$-sheaves, see \cite{lurie}, but to keep this note as far as possible at an informal level we will content ourselves with noticing that, for any open subset $U$ of $X$, there is a natural transformation of $\infty$-groupoids induced by the dgla map ${\mathcal T}_X(U) \xrightarrow{(\lie, e^{\boldsymbol{i}})}\holim \left(\xymatrix{ {\mathcal E}nd^{\geq 0}(\Omega_X^*)(U) \ar@<2pt>[r]^{{\rm incl.}}\ar@<-2pt>[r]_{0} & {\mathcal E}nd^{*}(\Omega_X^*)}(U)\right)$.}
\[ \Def({\mathcal T}_X) \xrightarrow{(\lie, e^{\boldsymbol{i}})}\holim \left(\xymatrix{ \Def({\mathcal E}nd^{\geq 0}(\Omega_X^*)) \ar@<2pt>[r]^{{\rm incl.}}\ar@<-2pt>[r]_{0} & \Def({\mathcal E}nd^{*}(\Omega_X^*))}\right).
\]
The homotopy fiber on the right should be thought as a homotopy flag manifold. Let us briefly explain this. At least na\"{\i}vely, the functor $\Def({\mathcal E}nd^{*}(\Omega_X^*))$ describes the infinitesimal deformations of the differential complex $\Omega_X^*$, whereas the functor $\Def({\mathcal E}nd^{\geq0}(\Omega_X^*))$ describes the deformations of the filtered complex $(\Omega_X^*,F^\bullet\Omega_X^*)$, i.e., of the pair consisting of the complex $\Omega_X^*$ \emph{and} the filtration $F^\bullet\Omega_X^*$. Therefore, the holimit describes a deformation of the pair (complex, filtration) together with a trivialization of the deformation of the complex. Summing up, the contraction of differential forms with vector fields induces a map of deformation functors
\[
\Def({\mathcal T}_X) \rightarrow {\rm hoFlag}(\Omega_X^*;F^\bullet\Omega_X^*),
\]
which we will call the \emph{local periods map} of $X$.

To recover from this the classical periods map, we just need to take global sections. Clearly, since we are working in homotopy categories, these will be derived global sections. The morphism of sheaves ${\boldsymbol{i}}: \T_X \to \Eps nd^*(\Omega_X^*)[-1]$ induces a Cartan homotopy ${\boldsymbol{i}}:{\bf R}\Gamma\T_X\to {\bf R}\Gamma\Eps nd^*(\Omega_X^*)[-1]$; composing this with the dgla morphism ${\bf R}\Gamma\Eps nd^*(\Omega_X^*)\to \End^*({\bf R}\Gamma\Omega_X^*)$ induced by the action of (derived) global sections of the endomorphism sheaf of $\Omega_X^*$ on (derived) global sections of $\Omega_X^*$, we get a Cartan homotopy
\[
{\boldsymbol{i}}:{\bf R}\Gamma\T_X\to \End^*({\bf R}\Gamma\Omega_X^*)[-1].
\]
The image of the corresponding dgla morphism $\lie$ (the derived globalization of Lie derivative) preserves the filtration $F^\bullet{\bf R}\Gamma\Omega_X^*$ induced by $F^\bullet\Omega_X^*$, so we have a natural map of $\infty$-groupoids
\[
\Def({\bf R}\Gamma\T_X)\to {\rm hoFlag}({\bf R}\Gamma\Omega_X^*;F^\bullet{\bf R}\Gamma\Omega_X^*)
\]
and, at the zeroth level, a map of {\bf Set}-valued deformation functors
\[
{\mathcal P}\colon\pi_{\leq 0}\Def({\bf R}\Gamma\T_X)\to \pi_{\leq 0}{\rm hoFlag}({\bf R}\Gamma\Omega_X^*;F^\bullet{\bf R}\Gamma\Omega_X^*)
\]
The functor on the left hand side is the {\bf Set}-valued functor of (classical) infinitesimal deformations of $X$; let us denote it by $\Def_X$. If we denote by $\End^*({\bf R}\Gamma\Omega_X^*;F^\bullet{\bf R}\Gamma\Omega_X^*)$ the subdgla of $\End^*({\bf R}\Gamma\Omega_X^*)$ consisting of endomorhisms preserving the filtration, then the pair $(\End^*({\bf R}\Gamma\Omega_X^*), \End^*({\bf R}\Gamma\Omega_X^*;F^\bullet{\bf R}\Gamma\Omega_X^*))$ is formal.\footnote{This is essentially a consequence of the $E_1$-degeneration of the Hodge-to-de Rham spectral sequence, see, e.g., 
\cite{deligne-illusie, faltings}.} Moreover, $
H^0(\End^*({\bf R}\Gamma\Omega_X^*))=\End^0(H^*_{dR}(X;{\mathbb K}))$
and $
H^0(\End^*({\bf R}\Gamma\Omega_X^*; F^\bullet{\bf R}\Gamma\Omega_X^*))=\End^0(H^*_{dR}(X;{\mathbb K}); F^\bullet H^*_{dR}(X;{\mathbb K}))$,\break
where $F^\bullet H^*_{dR}(X;{\mathbb K})$ is the Hodge filtration on the algebraic de Rham cohomology of $X$. By results described in the previous section, this means that the preimage of the distinguished point $0$ under the map
\[
 \pi_{\leq 0}{\rm hoFlag}({\bf R}\Gamma\Omega_X^*;F^\bullet{\bf R}\Gamma\Omega_X^*)\to
 \pi_{\leq 0} \Def({\mathcal E}nd^{\geq 0}(\Omega_X^*))
\]
is the quotient set
\[
 \frac{\exp(\End^0(H^*_{dR}(X;{\mathbb K})))}{\exp(\End^0(H^*_{dR}(X;{\mathbb K}); F^\bullet H^*_{dR}(X;{\mathbb K})))}\]
 and we recover the classical periods map of $X$
 \[
 {\mathcal P}\colon\Def_X\to{\rm Flag}(H^*_{dR}(X;{\mathbb K}); F^\bullet H^*_{dR}(X;{\mathbb K})).
\]
Also, the differential of ${\mathcal P}$ is the map induced in cohomology by the contraction of differential forms with vector fields,
\[
H^1({\boldsymbol{i}})\colon H^1(X,{\mathcal T_X})\to \int_p\Hom^0\left(F^pH^*_{dR}(X;{\mathbb K});\frac{H^*_{dR}(X;{\mathbb K})}{F^p H^*_{dR}(X;{\mathbb K})}\right),
\]
a result originally proved by Griffiths \cite{griffiths}. In the above formula, $\int_p$ denotes the end of the diagram
\[
\Hom^0\left(F^pH^*_{dR};\frac{H^*_{dR}}{F^p H^*_{dR}}\right)
\rightarrow
\Hom^0\left(F^pH^*_{dR};\frac{H^*_{dR}}{F^{p+1} H^*_{dR}}\right)
\leftarrow
\Hom^0\left(F^{p+1}H^*_{dR};\frac{H^*_{dR}}{F^{p+1} H^*_{dR}}\right)
\]
Also, we have the following version of the so-called Kodaira principle (ambient cohomology annihilates obstructions): obstructions to classical infinitesimal deformations of $X$ are contained in the kernel of
 \[
H^2({\boldsymbol{i}})\colon H^2(X,{\mathcal T_X})\to \int_p\Hom^1\left(F^pH^*_{dR}(X;{\mathbb K});\frac{H^*(X;{\mathbb K})}{F^p H^*_{dR}(X;{\mathbb K})}\right).
\]
In particular, if the canonical bundle of $X$ is trivial, then the contraction pairing
\[
H^2(X,{\mathcal T_X})\otimes H^{n-2}(X,\Omega^1_X)\to H^n(X;{\mathcal O}_X)\simeq {\mathbb K}
\] 
is nondegenerate, and so classical deformations of $X$ are unobstructed (Bogomolov-Tian-Todorov theorem, see \cite{bogomolov, tian, todorov}). Following \cite{iacono-manetti}, one immediately obtains the following refinement, due in its original formulation to Goldman and Millson \cite{goldman-millson}: if the canonical bundle of $X$ is trivial, then ${\bf R}\Gamma\T_X$ is a quasi-abelian dgla. To see this, just notice that the dgla map 
\[
{\bf R}\Gamma\T_X\xrightarrow{({\boldsymbol{l}},e^{\boldsymbol{i}})}
\holim \left(\xymatrix{ \End^*({\bf R}\Gamma\Omega_X^*;F^\bullet{\bf R}\Gamma\Omega_X^*) \ar@<2pt>[r]^{\phantom{mmm}{\rm incl.}}\ar@<-2pt>[r]_{\phantom{mmm}0} & \End^*({\bf R}\Gamma\Omega_X^*)}\right)
\]
is injective in cohomology and the target is a quasi-abelian dgla. Indeed, if $f\colon{\mathfrak g}\to{\mathfrak h}$ is a dgla morphism, with $H^*(f)$ injective and ${\mathfrak h}$ quasi-abelian, then the diagram of dglas
\par
\[
\xy (0,0)*{{\mathfrak g}\,} ; (15,-6)*+{{\mathfrak h}} **\dir{-}
?>* \dir{>}
,(30,0)*+{{\mathfrak k}} ; (15,-6)*+{{\mathfrak h}} **\dir{-}
?>* \dir{>}
,(30,0)*+{{\mathfrak k}} ; (45,-6)*+{V} **\dir{-}
?>* \dir{>}
,(24,-.6);(21,-3)**\dir{~}
,(39,-2.2);(36,-1.2)**\dir{~}
,(8,-1)*{\scriptstyle{f}}
 \endxy
\]
where $V$ is a graded vector space considered as a dgla with trivial differential and bracket, can be completed to a homotopy commutative diagram
\par
\[
\xy (0,0)*{{\mathfrak g}\,} ; (15,-6)*+{{\mathfrak h}} **\dir{-}
?>* \dir{>}
,(30,0)*+{{\mathfrak k}} ; (15,-6)*+{{\mathfrak h}} **\dir{-}
?>* \dir{>}
,(30,0)*+{{\mathfrak k}} ; (45,-6)*+{V} **\dir{-}
?>* \dir{>}
,(15,6)*+{{\mathfrak l}} ; (0,0)*+{{\mathfrak g}\,} **\dir{-}
?>* \dir{>}
,(15,6)*+{{\mathfrak l}} ; (30,0)*+{{\mathfrak k}} **\dir{-}
?>* \dir{>}
,(45,-6)*+{V} ; (60,-12)*+{W} **\dir{-}
?>* \dir{>}
,(24,-.6);(21,-3)**\dir{~}
,(39,-2.2);(36,-1.2)**\dir{~}
,(9,5.4);(6,3)**\dir{~}
,(8,-1)*{\scriptstyle{f}}
 \endxy
\]
with $W$ a graded vector space, and the composition ${\mathfrak l}\to W$ a quasi-isomorphism.
\par
As a conclusion, we recast the description of a period map for generalized deformations from \cite{fiorenza-manetti2} in the language of this note. Let $X$ be a smooth projective variety defined over the field $\mathbb{C}$ of the complex numbers, and denote by $\mathcal{P}oly^*_X$ the sheaf of dglas of \emph{multivector fields} on $X$, given by $
 \mathcal{P}oly^j_X=\bigwedge^{1-j}\mathcal{T}_X$, endowed with
 the zero differential and with the Schouten-Nijenhuis bracket.
Notice that $\mathcal{T}_X$ is a sub-Lie algebra of the dgla $\mathcal{P}oly^*_X$.
The contraction of differential forms with multivector fields 
\[
{\boldsymbol{i}}: \mathcal{P}oly^*_X \to \Eps nd^*(\Omega_X^*)[-1]
\]
 is a Cartan homotopy, and the corresponding 
dgla morphism $\lie$ is the Lie derivative along a multivector field, i.e., $\lie_\xi= [d_{\Omega^*_X}, {\boldsymbol{i}}_\xi]$. It is immediate that the image of $\lie$ is contained in the sub-sheaf of dglas:
\[ {\Eps nd_{0}^*(\Omega_X^*)} = \{f \in \Eps nd^*(\Omega_X^*)\, |\, f(\ker d_{\Omega^*_X}) \subseteq \mathrm{Im}(d_{\Omega^*_X}) \} \subset \Eps nd^*(\Omega_X^*),   \]
and so we have a natural transformation:
\[ \Def(\mathcal{P}oly^*_X) \xrightarrow{(\lie, e^{\boldsymbol{i}})}\holim \left(\xymatrix{ \Def({\Eps nd_{0}^*(\Omega_X^*)}) \ar@<2pt>[r]^{{\rm incl.}}\ar@<-2pt>[r]_{0} & \Def({\mathcal E}nd^{*}(\Omega_X^*))}\right),
\]
which we can think of as a local period map for generalized deformations. As above, to go from local to global, we take the derived global sections; then, taking $\pi_{\leq 0}$ we  obtain a natural morphism of {\bf Set}-valued deformation functors:
\[ \pi_{\leq 0}\Def(\mathbf{R}\Gamma\mathcal{P}oly^*_X) \xrightarrow{(\lie, e^{\boldsymbol{i}})}\pi_{\leq 0}\holim \left(\xymatrix{ \Def({\End_{0}^*({\bf R}\Gamma\Omega_X^*)}) \ar@<2pt>[r]^{{\rm incl.}}\ar@<-2pt>[r]_{0} & \Def(\End^{*}({\bf R}\Gamma \Omega_X^*))}\right).
\]  
On the left,  $\pi_{\leq0} \Def(\mathbf{R}\Gamma\mathcal{P}oly^*_X)$ is the functor $\widetilde{\Def}_X$  of generalized deformations of $X$. It is shown in \cite{fiorenza-manetti2}, using the Dolbeault resolution as a model for $\mathbf{R}\Gamma\Omega^*_X$, and making use of the $\partial\overline{\partial}$-lemma, that 
the pair $(\End^*_0(\mathbf{R}\Gamma\Omega^*_X),\End^*(\mathbf{R}\Gamma\Omega^*_X))$ is quasi-isomorphic to the pair $(0,\End^*(H^*(X,\mathbb{C})))$. Hence, one obtains the period map for generalized deformations:
\[
\widetilde{\mathcal{P}}:\widetilde{\Def}_X\to \exp(\End^0(H^*(X,\mathbb{C})).
\]
The tangent map $d\widetilde{\mathcal{P}}$ is the contraction of differential forms with multivector fields, read at the cohomology level:
\[
H^1(\boldsymbol{i}):\bigl(\oplus_kH^k(X;\wedge^k\mathcal{T}_X)\bigr)\otimes \bigl(\oplus_{p,q}H^q(X,\Omega^p_X)\bigr)\to \oplus_{p,q,k} H^{q+k}(X,\Omega^{p-k}_X),
\]
and obstructions to generalized deformations are contained in the kernel of the contraction
\[
H^2(\boldsymbol{i}):\bigl(\oplus_kH^{k+1}(X;\wedge^k\mathcal{T}_X)\bigr)\otimes \bigl(\oplus_{p,q}H^q(X,\Omega^p_X)\bigr)\to \oplus_{p,q,k} H^{q+k+1}(X,\Omega^{p-k}_X).
\]
In particular, from this one recovers Barannikov-Kontsevich's result, that generalized deformations of a smooth projective Calabi-Yau manifold are unobstructed \cite{barannikov-kontsevich}. 
\par
It is tempting to extend the construction of the period map for generalized deformations to the case of a smooth projective manifold defined on an arbitrary characteristic zero algebraically closed field $\mathbb{K}$, 
\[
\widetilde{\mathcal{P}}:\widetilde{\Def}_X\to \exp(\End^0(H^*_{dR}(X;\mathbb{K})).
\]
To do this one only has to prove that $(\End^*_0(\mathbf{R}\Gamma\Omega^*_X),\End^*(\mathbf{R}\Gamma\Omega^*_X))$ is quasi-isomorphic to $(0, \End^*(H^*_{dR}(X,\mathbb{K}))$. Yet, to mimic the argument in \cite{fiorenza-manetti2} one needs an algebraic substitute of the $\partial\overline{\partial}$-lemma. A natural candidate for this is $E_1$-degeneracy of the Hodge-to-de Rham spectral sequence for a smooth projective manifold. It has however to be remarked that, while in the
$\partial\overline{\partial}$-lemma the two differentials $\partial$ and
$\overline\partial$ play perfectly interchangeable roles, see, e.g. \cite[Corollary 3.2.10]{huybrechts}, this is not true for the \v{C}ech
and the de Rham differentials in the \v{C}ech-de Rham bicomplex $\check{C}^q(\mathcal{U},\Omega^p_X)$ associated with an open cover $\mathcal{U}$ of $X$. In particular only one of the two spectral sequences associated with this bicomplex, namely the 
Hodge-to-de Rham spectral sequence, degenerates at $E_1$. This asymmetry seems to suggest that a purely algebraic proof of the quasi-isomorphism $(\End^*_0(\mathbf{R}\Gamma\Omega^*_X),\End^*(\mathbf{R}\Gamma\Omega^*_X))\simeq (0, \End^*(H^*_{dR}(X,\mathbb{K}))$ could be a nontrivial result.

\end{document}